\documentclass[letterpaper,12pt]{article}

\title{Computing the Integer Programming Gap}

\author{Serkan Ho{\c s}ten  \\
{\small Department of Mathematics, 
San Francisco State University, San Francisco} \\
\and Bernd Sturmfels \!\!\!\! \thanks{Partially supported by the
National Science Foundation (DMS-0200729)}
 \\
{\small Department of Mathematics, University of California,
Berkeley} }

\date{}

\usepackage{latexsym,array,delarray,amsthm,amssymb,psfrag}
\usepackage[sumlimits]{amsmath}
\usepackage[dvips]{graphicx}

\theoremstyle{plain}
\newtheorem{thm}{Theorem}[section]
\newtheorem{lemma}[thm]{Lemma}
\newtheorem{prop}[thm]{Proposition}
\newtheorem{cor}[thm]{Corollary}

\theoremstyle{definition}

\newtheorem{ex}[thm]{Example}

\theoremstyle{remark}
\newtheorem*{rmk}{Remark}

\newcommand{\zz}{\mathbb{Z}}
\newcommand{\nn}{\mathbb{N}}

\newcommand{\qq}{\mathbb{Q}}
\newcommand{\rr}{\mathbb{R}}

\newcommand{\LL}{\mathcal{L}}

\begin{document}
\maketitle

\begin{abstract}
\noindent
We determine the maximal gap between the optimal values
of an integer program and its linear
programming relaxation, where the matrix and
cost function are fixed but the right hand side
is unspecified. Our formula involves irreducible
decomposition of monomial ideals. The gap can
be computed in polynomial time when the dimension is fixed.

\end{abstract}

\section{Introduction}

We consider the general integer programming problem in standard form,
\begin{equation}
\label{ip}
{\rm Minimize} \quad c \cdot u \quad
\hbox{subject to} \quad A \cdot u \, = \, b, \, u \geq 0, \, u \,\, \hbox{integral}.
\end{equation}
Here $A$ is a fixed $d \times n$ integer matrix, $b \in \zz^d$ and
$c \in \qq^n$. Its linear programming relaxation is obtained by
dropping the integrality constraints,
\begin{equation}
\label{lp}
{\rm Minimize} \quad c \cdot u \quad
\hbox{subject to} \quad A \cdot u \, = \, b, \, u \geq 0.
\end{equation}
Suppose that the integer program (\ref{ip}) is feasible and
bounded. Then the linear program (\ref{lp}) is feasible
and bounded as well, and the optimal value of (\ref{ip}) is greater than
or equal to the optimal value of (\ref{lp}). We define
the nonnegative rational number $\,{\rm gap}(A,c) \,$
to be the maximum difference of the two optimal values
as $b$ ranges over vectors such that $(\ref{ip})$
is feasible and bounded.  It follows from known results
\cite[Theorem 17.2]{Sch} that this maximum is
bounded.

Our main aim in this paper is
to provide an exact formula for $\,{\rm gap}(A,c) $.
We express our results using the language of
Gr\"obner bases, as in \cite[Chapter 8]{CLO}, \cite{HT1},
\cite[\S 4.4]{SST}.
A nonnegative integer vector $\,u = (u_1,\ldots,u_n) \in \nn^n\,$
is called {\it non-optimal} if it is not an optimal solution
of (\ref{ip}) with $\, b = A \cdot u $.
We represent each non-optimal vector $u$ by
a monomial $\,x^u =  x_1^{u_1} x_2^{u_2} \cdots x_n^{u_n}$,
and we consider the ideal $\,M(A,c) \,$ generated
by these monomials in the polynomial ring $\, \rr[x_1,\ldots,x_n]$.
The minimal generators of $\,M(A,c)\,$ can be read off
from a   Gr\"obner basis for (\ref{ip}). If $c$ is generic
then $\, M(A,c) \,$ is an initial ideal
of the toric ideal of $A$; see \cite{GB+CP}.
If $c$ is not generic then we can compute
$\, M(A,c) \,$  using \cite[Algorithm 4.4.2]{SST}.
 A monomial ideal $I$
in $\rr[x_1,\ldots,x_n]$ is called {\it irreducible} if it is
generated by powers of the variables:
$$ I  \quad = \quad  \langle \,
x_{i_1}^{u_{i_1}+1},\, x_{i_2}^{u_{i_2}+1},\, \ldots , \,x_{i_r}^{u_{i_r}+1}\, \rangle .$$
Every monomial ideal $M$  in $\rr[x_1,\ldots,x_n]$
  can be written uniquely as an irredundant
intersection of finitely many irreducible monomial ideals, which are called
the {\it irreducible components} of $M$.
Suppose that $\, I \,$ is an irreducible component of
$M(A,c)$.
We define the {\it gap value}  of $I$
with respect to $A$ and $c$ to be the optimal objective function
value of the auxiliary linear program
\begin{equation}
\label{gapprog}
\begin{array}{ll}
{\rm Maximize} & u_{i_1} c_{i_1} \, + \,
   u_{i_2} c_{i_2} \, + \cdots + \, u_{i_r} c_{i_r} \,\, - \,\, c \cdot v   \\
\hbox{subject to} & A \cdot v \quad = \quad  u_{i_1} \cdot A_{i_1}
\, +\,  u_{i_2}  \cdot A_{i_2} \,+ \cdots +\,  u_{i_r} \cdot A_{i_r},   \\
 & \hbox{and} \quad  v_{i_1}, \ldots, v_{i_r} \geq 0.
\end{array}
\end{equation}
Here $A_1, A_2,\ldots,A_n$ are the column vectors of the matrix $A$.

\begin{thm} \label{firstthm}
The integer programming gap, $ {\rm gap}(A,c) $,
equals the maximum gap value of any
irreducible component $\, I \,$ of the monomial
ideal $M(A,c)$.
\end{thm}

We remark that $\, {\rm gap}(A,c) \,$ is zero if and
only if the monomial ideal $M(A,c)$ is
generated by squarefree monomials 
$\, x_{j_1} x_{j_2} \cdots x_{j_r} $.

This paper is organized as follows. In Section 2
we rephrase our problem in the more general
setting of  lattice programs, and we
prove Theorem \ref{firstthm} in this context.
In Section 3 we apply the work
of Barvinok and Woods \cite{BW}
on short rational generating functions to derive the
following complexity result.

\begin{thm} \label{complexity}
For fixed $d$ and $n$, the
integer programming gap, $ {\rm gap}(A,c) $,
can be computed in polynomial time
in the binary encoding of $A$ and $c$.
\end{thm}

Section 4 concerns applications to the
statistical theory of multidimensional
contingency tables. Here we are interested
in the integer programming gap of certain higher-dimensional
transportation problems. These play an important role in data security.
For the statistical background see \cite{CG}, \cite{DF}.

In Section 5 we vary the cost vector $c$,
and we prove that the function
\begin{equation}
\label{gapfct}
{\rm gap}_A \,\, :\,\, \rr^n \rightarrow \rr,\,
c \mapsto {\rm gap}(A,c)
\end{equation}
is a piecewise-linear function
on $\rr^n$. We show that ${\rm gap}_A$ is linear on the
cones of a fan which refines the familiar \emph{Gr\"obner fan}
(cf.~\cite{HubT} and \cite{ST}).

\smallskip

We close the introduction  with a small example.
Let (\ref{ip}) be the problem of \emph{making change}
using  pennies, nickels, dimes and quarters,
where the number of
coins is fixed, and nickels and quarters are
used most sparingly. In symbols,
$$ d = 2 \, , \,\, n= 4, \quad
A =
\left[ \begin{array}{cccc}
1 & 1 & 1 & 1 \\
1 & 5 & 10 & 25
\end{array} \right], \quad
c =
\left[ \begin{array}{cccc}
 0 & 1 & 0 & 1
\end{array} \right]. $$
Using notation as in \cite[\S I.5]{EGSS},
this problem is solved by the Gr\"obner basis
$$ {\cal G} \quad = \quad
\bigl\{ \,\,
\underline{n^3 q} \,-\, d^4 \,, \,\,
\underline{n^6} \,- \,p^5 q \,, \,\,
\, \underline{n^3 d^4} \,- \,p^5 q^2 \,, \,\,
\underline{p^5 q^3}-d^8 \, \bigr\}. $$
The four underlined leading monomials generate the ideal
\begin{equation}
\label{earliermono}
 M(A,c) \quad = \quad
\langle \, n^3,p^5 \,\rangle \,\, \cap \,\,
\langle \, n^3, q^3 \,\rangle  \,\, \cap \,\,
\langle \, n^6, d^4, q \,\rangle.
\end{equation}
The gap values of the three irreducible components are
$\, 76/15, \, 4 $ and $ 5 $. Hence
$$ {\rm gap}(A,c) \quad = \quad  76/15 \quad = \quad 5.0666666... $$
This gap is attained when
expressing one dollar and $14$ cents with ten coins.
The optimal solutions are $(4,2,0,4)$
and $(0,0, 136/15, 14/15)$ respectively.

\section{The gap theorem for lattice programs}

Lattice programs are defined as follows.
Let  $\LL $ be a fixed lattice of rank $m$ in $\zz^n$.
Then  $\,\LL_{\rr} =  \LL \otimes_{\zz} \rr \,$ is an
$m$-dimensional  vector space in $\rr^n$.
We also fix a cost vector $c \in \qq^n$.
 Now, for every $z \in \nn^n$, we get a
{\em lattice program}
\begin{equation}
\label{latprog}
{\rm Minimize} \quad  c \cdot v \quad
\hbox{subject to} \quad v \equiv z \quad \hbox{mod } \, \LL,
\quad v \in \nn^n ,
\end{equation}
and its {\em linear relaxation}
\begin{equation}
\label{latproglp}
{\rm Minimize} \quad  c \cdot v \quad
\hbox{subject to} \quad v \equiv z \quad \hbox{mod } \, \LL_{\rr},
\quad v \geq 0   .
\end{equation}

We define the \emph{lattice
programming gap}  $\, {\rm gap}(\LL,c) \,$
to be maximum of the differences between the optimal values
of (\ref{latprog}) and (\ref{latproglp}) as $z$ runs over $\nn^n$.

\begin{rmk} The integer programs (\ref{ip}) are lattice programs
with $\LL = \ker(A) \cap \zz^n$ and $b = A \cdot z$. The linear
programming relaxations (\ref{lp}) correspond to  the
linear relaxations (\ref{latproglp}) of these lattice programs.
Note that we have $m = n -d$.
\end{rmk}

Let $M(\LL, c)$ be the ideal
generated by all monomials $\, x^u \, $
where $u\in \nn^n$ is non-optimal for
(\ref{latprog}).  Any irreducible component of
$M(\LL, c)$ has the form
$$ I(u,\tau) \quad = \quad \langle \,
x_{i_1}^{u_{i_1}+1},\, x_{i_2}^{u_{i_2}+1},\, \ldots , \,x_{i_r}^{u_{i_r}+1}\,
\rangle, $$
where  $\tau = \{i_1, \ldots, i_r\}$ and $u \in \nn^n$ with $u_j = 0$ for
$j \notin \tau$.
The \emph{gap value} of $I(u,\tau)$ with respect to
$\LL$ and $c$ is the optimal objective function value of
\begin{equation}
\label{latgapprog}
\begin{array}{ll}
{\rm Maximize} & u_{i_1} c_{i_1} \, + \, u_{i_2} c_{i_2} \, + \cdots + \, 
u_{i_r} c_{i_r} \,\, - \,\, c \cdot v   \\
\hbox{subject to} \quad & v \equiv u \quad \hbox{mod } \, \LL_{\rr}
\quad \hbox{and} \quad  v_{i_1}, \ldots, v_{i_r} \geq 0
\end{array}
\end{equation}
We restate  Theorem \ref{firstthm} for the more general
setting of lattice programs.

\begin{thm}
\label{latmainthm}
The lattice programming gap, ${\rm gap}(\LL,c) $,
equals the maximum  gap value of any
irreducible component $ I(u,\tau) $ of the monomial
ideal $ M(\LL,c)$.
\end{thm}

\begin{ex}
\label{finiteindex} $(m=n)$ \
Let ${\cal L}$ be a finite index sublattice of $\zz^n\,$
and $c$ a nonnegative vector.  Since $\LL_{\rr} = \rr^n$,
the objective function value of (\ref{latproglp})
is always zero, so  $\, {\rm gap}(\LL,c) \,$
is the largest objective function value of (\ref{latprog}).
The finite abelian group $\zz^n/\LL$ is in bijection with
the set of monomials not in $M(\LL,c)$. The
irreducible components $\,I(u,\{1,\ldots, n\}) \,$  are
indexed by monomials $\,x_1^{u_1} x_2^{u_2} \cdots x_n^{u_n}\,$
which are maximal with respect to divisibility
among those not in $M(\LL,c)$.
The lattice programming gap
is the maximum of the corresponding values
$\,c_1 u_1 +  c_2 u_2 + \cdots +  c_n u_n $. \qed
\end{ex}

Example  \ref{finiteindex} is called the
\emph{group problem} in the integer programming
literature \cite[\S 24.2]{Sch}.
More generally, every lattice program  (\ref{latprog})
has a natural family of  {\em group relaxations} which
are indexed by subsets $\tau$ of $\{1,2,\ldots,n\}$:
\begin{equation}
\label{group}
{\rm Minimize} \quad  c \cdot v \quad
\hbox{subject to} \quad v \equiv z \quad \hbox{mod } \, \LL,
\quad v \in \zz^n, \,\, v_i \geq 0  \,\,\, \forall \, i \in \tau.
\end{equation}

If an optimal solution $v^*$ of (\ref{group}) is a nonnegative vector
then $v^*$ is also feasible and optimal for (\ref{latprog}).
In this case we say that lattice program (\ref{latprog})
is \emph{solved by} $\tau$. The minimal collection
of required subsets $\tau$ is studied in
\cite{HT2}. The following result is well-known in
the algebraic theory of integer programming;
see \cite[\S 3]{HT1}. The pairs
$(u,\tau)$ in Proposition \ref{stdpairs}
are called \emph{standard pairs}.
A combinatorial introduction can be found in \cite{T}.  

\begin{prop}
\label{stdpairs}
There is a unique minimal finite set $\mathcal{S}$ of
irreducible ideals $I(u,\tau)$ whose intersection is $M(\LL, c)$
such that every optimal solution $v^*$ to (\ref{latprog}) has the form
 $v^* = u + v'$,  with $v_i' = 0$ for $i \in \tau$, for some
$I(u, \tau) \in \mathcal{S}$.
\end{prop}

When a lattice program  (\ref{latprog})
has an optimal solution $v^* = u + v'$ as in Proposition \ref{stdpairs},
we shall say that it is \emph{solved by the standard ideal} $I(u, \tau)$.
In this case the group relaxation (\ref{group}) has the same
optimal solution $v^*$.

\begin{rmk} The irreducible components of $M(\LL, c)$ 
are a subset of $\mathcal{S}$. In fact,
the irreducible components are the minimal elements of $\mathcal{S}$ when
the standard ideals  are ordered with respect to inclusion.
In the special case of Example  \ref{finiteindex}, we have
$\mathcal{S} \, = \, \bigl\{ I(u,\, \{1,\ldots, n\})\, : \, x^u \notin
M(\LL,c) \bigr\}$.
\end{rmk}

\begin{lemma}
\label{lem1} Fix $I(u, \tau) \in \mathcal{S}$.
The gap value of $I(u, \tau)$ equals
the maximum difference between the optimal values of
 (\ref{latprog}) and (\ref{latproglp})
as $z$ ranges over all vectors in $\nn^n$ such that the program 
(\ref{latprog}) is solved by the standard ideal $I(u,\tau)$.
\end{lemma}
\begin{proof}
Suppose a lattice program (\ref{latprog}) is solved by
the standard ideal $I(u, \tau)$  and has
the optimal solution $x^* = u + x'$ where $x_i' = 0$ for all
$i \in \tau$.
Let $y^*$ be an optimal solution for the linear relaxation
(\ref{latproglp}). Then the  difference between the optimal  values is
$c \cdot x^* - c \cdot y^*$.
Since $y^* - x'$ is a feasible solution for (\ref{latgapprog}),
the optimal value of this program is at least $c \cdot u -
c \cdot (y^* - x') = c \cdot x^* - c \cdot y^*$.

Hence we only need  to find a vector $z \in \nn^n$ whose associated
lattice program (\ref{latprog})  is solved by $I(u, \tau)$, and such that
the difference of the optimal values of  (\ref{latprog}) and (\ref{latproglp})
is greater than or equal to the optimal value  of (\ref{latgapprog}).
 Let $v^*$ be an optimal solution to
(\ref{latgapprog}) and define a vector $v' \in \nn^n$
by $\, v'_i \, = \, {\rm max} \bigl\{ 0, -\lfloor v_i^* \rfloor \bigr\} $.
Then the  lattice program (\ref{latprog}) with $z = u + v'$
is solved by $I(u, \tau)$. In fact,  $z$ is the optimal
solution, and $v^* + v'$ is a feasible
solution for the linear relaxation (\ref{latproglp}).
The difference between the optimal solution values
of the two programs is at least $ c \cdot (u + v') - c \cdot (v^* + v') =
c \cdot u - c \cdot v^*$.
\end{proof}

\noindent
{\sl Proof of Theorem \ref{latmainthm}}: In light of Lemma \ref{lem1}
and Proposition \ref{stdpairs}, we just need to show that
if $I(u, \tau)$ and $I(u', \tau)$ are two standard ideals
with $u \leq u'$, then the optimal value of (\ref{latgapprog})
is at most that of (\ref{latgapprog}) with $u$ replaced by $u'$.
In order
to do this we will reformulate these programs.
Let $B$ be an $n \times m$ matrix
whose columns form a lattice basis of $\LL$, and let
$\{b_1, \ldots, b_n\} \subset \zz^m$ be the rows of $B$.
The feasible solutions to (\ref{latgapprog}) are in bijection
with $\{t \in \rr^m \, : \, b_i \cdot t \leq u_i, \,\, i \in \tau \}$
via $t \mapsto v = u - Bt$. If we let $w := \sum_{i=1}^n c_ib_i$,
then the following linear program is equivalent to (\ref{latgapprog})
and has the same objective function value:
\begin{equation}
\label{newgap}
{\rm Maximize} \quad  w \cdot t \quad
\hbox{subject to} \quad b_i \cdot t \leq u_i \quad \forall i \in \tau
\end{equation}
The cost vector $w $ is independent of $u$.
If we replace $u$ by $u'$ in (\ref{newgap}) then  the
feasible region increases but the objective function is unchanged.
Hence the optimal value of 
(\ref{newgap}) can only increase when replacing $u$ by $u'$.
 \hfill $\Box$

\begin{ex} Theorem \ref{latmainthm} suggests that we compute
${\rm gap}(\LL, c)$ by solving finitely many
linear programs (\ref{latgapprog}),
one for each irreducible component $I(u, \tau)$ of $M(\LL, c)$.
One difficulty is that the number of irreducible components
can be exponential in the problem size.
We  illustrate this phenomenon with an example
taken from \cite[\S 4]{SWZ}. For $r \geq 4$ let
$\LL_r$ be the lattice  generated by
$$(r, \, r, \, r), \,\,\,  (r-1, \, r+1, \, r-1) \,\,\,
\hbox{and}\, \,\, (0, \, 0, \, r-2) \quad \hbox{in} \,\,\,
\zz^3. $$
The index of $\LL_r$ in $\zz^3$ is  $2r(r-2)$, so we
 are in the situation of Example \ref{finiteindex}.
Let $c$ be a cost vector that chooses the
{\em degree lexicographically} smallest solution of (\ref{latprog}).
{}From the Gr\"obner basis in  \cite[Lemma 4.5]{SWZ} we find that
$$ M(\LL_r, c) \quad = \, \,
\mathop{\bigcap_{1 \leq a \leq r-2}}_{b = 2r+1-a}
\langle x^a, \, y^b, \, z \rangle   \,\,\,\,  \cap \,\,\,\,
\mathop{\bigcap_{1 \leq a \leq r-3}}_{b = r-1-a}
\langle x, \, y^a, \, z^b \rangle. $$
The number of irreducible components is $2r-5$. This
is exponential in $\, O({\rm log}(r)) $, the bit complexity
of the data. In the next section we demonstrate how
the monomial ideals $ M(\LL, c) $ can be encoded more
efficiently.
\qed
\end{ex}

\section{Gap function as a short rational function}

In this section we prove Theorem \ref{complexity}. The result
extends easily to the lattice programs of Section 2,
but for the sake of notational convenience, we present
the relevant generating functions in the  original setting
of integer programs (\ref{ip}). Throughout the section
we assume that $d$ and $n$ are fixed integers.

Let $\nn A$ denote the semigroup in $\zz^d$ spanned by
the columns $A_1, \ldots, A_n$ of the matrix $A$.
The elements of $\nn A$ are the feasible right-hand-side vectors
of the integer programs we consider. Let
$\mu_{IP} : \nn A \mapsto \qq$ be the function whose
value at $b$ is the optimal value of
the integer program (\ref{ip}), and let
$\mu_{LP} : \nn A \mapsto \qq$ be the corresponding
optimal value function for the linear program (\ref{lp}).
Since $c$ is assumed to have rational coordinates, we
can compute (in polynomial time)
a \emph{global denominator} $\Delta \in \nn$
by multiplying the  least common multiple
of the maximal minors of $A$ with the
least common denominator of $c_1,c_2,\ldots,c_n$.
This choice of $\Delta$ ensures that every value of
the functions $\mu_{IP}$ or $\mu_{LP}$
is an integer multiple of $\Delta$.

We are interested in the  generating function
$$ G(t_1, \ldots, t_d; s) \,\,\, := \, \,\,
\sum_{b \in \nn A} {\bf t}^bs^{\mu_{IP}(b) - \mu_{LP}(b)}
\quad \in \,\,\, \qq[s^{1/\Delta}][[ \nn A]]. $$
The ambient ring is the completion of the
semigroup algebra of $\nn A$ with coefficients
in the univariate polynomial ring $\qq[s^{1/\Delta}]$.
We shall see that $G$ is a rational series
which can be represented by an element of
$$ \qq(t_1,\ldots,t_d)[s^{1/\Delta}] \quad \subset \quad 
\qq(t_1, \ldots, t_d, s^{1/\Delta}).  $$
Thus $G$ is a polynomial in $s^{1/\Delta}$ whose coefficients
are rational functions in $t_1,\ldots,t_d$.
The degree of that univariate polynomial is the
integer programming gap $\, {\rm gap}(A,c) $.
We shall prove the following complexity result.

\begin{thm}
\label{genfct}
The rational function
$G(t_1, \ldots, t_d; s) $
can be computed in polynomial time
in the binary encoding of the matrix $A$
and the cost vector $c$.
\end{thm}

What we are claiming is that $G({\bf t};s)$ is a
\emph{short rational generating function}
in the sense of Barvinok and Woods \cite{BW},
and our proof is a direct application of their work.
A different approach to the integer programming gap using
generating functions and integer programming duality is presented
by Lasserre in \cite{L}.
We illustrate the main point of Theorem \ref{genfct} in an example.

\begin{ex}
Let $d = 1, n =2, a,b \in \nn$ and consider the integer program
$$ \hbox{ minimize } \,\,
u_1 \,\, \hbox{ subject to }\,\,  u_1 + a u_2  = b, \,\,
u_1, u_2 \geq 0.  $$
The optimal values are
 $\mu_{LP}(b) = 0$ and
 $\mu_{IP}(b) = b \,\,{\rm mod}\,\, a$. Hence
$$ G(t,s) \quad  = \quad
\sum_{b \in \nn} t^b s^{b \,{\rm mod} \, a}
\quad = \quad
\sum_{i=0}^{a-1}
\sum_{j=0}^\infty t^{i+ja} s^i \\
 = \quad
\sum_{i=0}^{a-1}
(t^i/(1-t^a)) \cdot s^i. $$
This is a polynomial in $s$ with $a$ terms. Its expansion
requires exponential space in the bit size
$\, O({\rm log}(a))\, $
of the given data $\, A = [1 \,\, a] \,$ and $c = [1 \,\, 0 ]$.
On the other hand, clearly this gap polynomial
is a short rational function
$$ G(t,s) \quad  = \quad
\frac{1 - s^a t^a}{ (1-t^a)(1-st)}. $$
Our strategy is to compute $G(t,s)$ in polynomial time
and then to extract
$\, {\rm gap}(A,c) \, = \,
{\rm degree}_s \bigl( G(t,s) \bigr) \, = \,  a-1$. \qed
\end{ex}

One ingredient in the proof of Theorem \ref{genfct}
is the {\em Hadamard product} $*$ of two generating
functions: if $g_1({\bf x}) = \sum_u \beta_u {\bf x}^u$
and $g_2({\bf x}) = \sum_v \gamma_v {\bf x}^v$ then
$g_1({\bf x}) * g_2({\bf x})$ is the generating function
$\,\sum \beta_u \gamma_u {\bf x}^u$. The proof of
\cite[Lemma 3.4]{BW} tells us that, if $g_1({\bf x})$ and
$g_2({\bf x})$ are short rational functions, then their
Hadamard product $g_1({\bf x}) * g_2({\bf x})$
can be computed in polynomial time.  Another
ingredient is the following lemma which
is of independent interest.

\begin{lemma} \label{shortopt}
The generating function for all the optimal points,
$$ H({\bf x}) \quad = \quad \sum \bigl\{ \,{\bf x}^u
\,\, : \,\, u \, \hbox{is optimal for  (\ref{ip}) with} \,\,
b = Au \, \bigr\}. $$
can be computed in polynomial time,
in the binary encoding of $A$ and $c$.
\end{lemma}

\begin{proof}
 The proof is an adaptation
of \cite[\S 7.3]{BW}. Without loss of generality
we assume that $c \in \zz^n$. Let $S = \{u: \, x^u \in M(A, c)\}$
be the set of non-optimal points,
and let $f(S;x_1, \ldots, x_n) = \sum_{u \in S} {\bf x}^u$ be the
generating function of $S$.  This generating function is a
rational function; in particular, we have
$$ f(S; {\bf x}) \prod_{i=1}^n (1 - x_i)
\quad = \quad g(S; {\bf x})$$
where $g(S; {\bf x})$ is a polynomial.  In view of the identity
$$ H({\bf x}) \quad = \quad
\frac{1 - g(S; {\bf x})}{\prod_{i=1}^n  (1 - x_i)} , $$
it suffices to show that $g(S; {\bf x})$
can be computed in polynomial time.

 We let $L := (n-d)D(A)$
where $D(A)$ is the maximum of the absolute values of the maximal
minors of $A$. Since we fix $d$ and $n$, the bit size of $L$
is a polynomial in the bit size of the data. Theorem 4.7 of
\cite{GB+CP} implies that if ${\bf x}^m$ is a term in $g(S; \bf{x})$
then $m_i \leq L$ for $i=1,\ldots, n$. We let
$\Lambda := \{(u_1, \ldots, u_n) \in \nn^n \, : \, u_i \leq L \,, \,
\forall \, i\}$ and $f(\Lambda; {\bf x})$ its generating function.
Furthermore let
$$ C \,\,\,:= \,\,\,
 \{(u,v) \in \nn^{2n} \, : \, Au = Av, \,\, cu \geq cv + 1, \,\,
\mbox{ and } u_i \leq L, \,\, v_i \leq 2L, \,\, \forall \, i\}$$
The projection of $C$ onto the first $n$ coordinates will
be denoted by $S'$, and we claim that $S' = S \cap \Lambda$. The
inclusion $S' \subseteq S \cap \Lambda$ is clear. For the other inclusion
let $u \in S \cap \Lambda$. The theory of Gr\"obner bases of toric
ideals \cite{GB+CP} implies that there exists $(u',v') \in \nn^{2n}$
with $u' \leq u$, $Au' = Av'$ and $cu' \geq cv'+ 1$, and
where $u_i' v_i' = 0$,
$u_i' \leq L$ and $v_i' \leq L$  for
$i =1, \ldots, n$. Now we let $v = u - u' + v'$ and observe
that $(u,v) \in C$. Since $S'$ is the projection of all lattice
points in a polytope, Theorem 1.7 in \cite{BW} implies that
$\, f(S'; {\bf x}) = \sum_{u \in S'} {\bf x}^u$, and hence
$\,g(S'; {\bf x}) = f(S'; {\bf x}) \prod_{i=1}^n (1-x_i)$, can be computed
in polynomial time. The claim we proved above says that
$$ f(S; {\bf x})  \quad = \quad f(S'; {\bf x})  \, + \,
\sum_{u \in S \backslash \Lambda} {\bf x}^u, $$
and this implies
$ \,g(S; {\bf x}) \, = \, g(S'; {\bf x}) \, + \, h({\bf x}) $,
where $h({\bf x})$ is a series none of whose terms has its
exponent vector in  $\Lambda$.
Now we conclude that
$$ g(S; {\bf x}) \quad = \quad g(S; {\bf x}) * f(\Lambda; {\bf x})
\quad = \quad
g(S'; {\bf x}) * f(\Lambda; {\bf x}).$$
The Hadamard product on the right can be
computed in polynomial time. \end{proof}

\noindent
{\sl Proof of Theorem \ref{genfct}:}
We replace each variable $x_i$ by ${\bf t}^{A_i}s^{c_i}$
in the rational function $ H({\bf x})$. This monomial
substitution  can be done in  polynomial time \cite[Theorem 2.6]{BW}.
The result is the short rational generating  function
$$ H_{IP}(t_1, \ldots, t_d; s) \, = \, \sum_{b \in \nn A} {\bf t}^b
s^{\mu_{IP}(b)}.$$

\noindent
The last series that is left to compute is
$$H_{LP}(t_1, \ldots, t_d; s) \, = \, \sum_{b \in \nn A} {\bf t}^b
s^{-\mu_{LP}(b)},$$
since $G({\bf t}; s) \, = \, H_{IP}( {\bf t}; s) *_{\bf t}
H_{LP}({\bf t}; s)$.
Let $\sigma = \{i_1, \ldots, i_d\} \subset \{1, \ldots, n\}$
be an optimal basis of (\ref{lp}) for some   $b$ in  $\nn A$.
The number of optimal bases $\sigma$ is constant since
$n$ and $d$ are fixed. The union
of $S_{\sigma} := \nn A  \cap \rr_{\geq 0}\{A_i \, : \, i \in \sigma\}$
as $\sigma$ varies over all optimal bases is equal to $\nn A$.
The generating function $f(S_\sigma; {\bf t}) \,  = \,
\sum_{b \in S_\sigma} {\bf t}^b$ can be computed in polynomial
time since $S_\sigma$ is the set of lattice points in a rational
polyhedral cone \cite{BP}. We let ${\hat c} \in \rr^d$ be
the unique vector such that $ (A^t \cdot {\hat c})_i = c_i$
for $i \in \sigma$. Now the generating function
$g_\sigma({\bf t}; s) \, = \, \sum_{b \in S_\sigma}
{\bf t}^b s^{-\mu_{LP}(b)}$ is obtained from $f(S_\sigma; {\bf t})$
by the monomial substitution $t_i = t_i s^{-{\hat c}_i}$ for
$i=1, \ldots, d$. Finally, we use \cite[Corollary 3.7]{BW}
to compute $H_{LP}({\bf t}; s)$ by patching the
series $g_\sigma({\bf t}; s)$ for the various  bases $\sigma$.
\hfill $\Box$

\smallskip
We now show that ${\rm gap}(A,c)$, which is
the degree of $G({\bf t};s)$ as a polynomial in $s$, can
be extracted in polynomial time. This uses the
following lemma.

\begin{lemma} \label{finddegree}
Let $f({\bf t}; s) \in \qq(t_1,\ldots,t_d)[s]$ be
a short non-zero rational function and $K$
a known upper bound on ${\rm deg}_s(f({\bf t}; s))$. Then
this degree can be computed with
$\,{\rm log}(K)\,$ Hadamard products of  $f({\bf t}; s)$ with
polynomials in $s$.
\end{lemma}

\begin{proof}
Without loss of generality we assume $K = 2^k$ for some $k \in \nn$. Let
$$g_{[p,r]}(s) \, = \, \sum_{i=p}^r s^i \, = \, \frac{s^p - s^{r+1}}{1-s},$$
and use the following binary search algorithm:

\noindent
{\tt findDegree}($f({\bf t}; s), p, r$)
\begin{itemize}
 \item[0.] If $p = r$ return $p$.
 \item[1.] If $f({\bf t}; s) * g_{[p,r]}(s)$ is not identically zero
then \\
{\tt findDegree}($f({\bf t}; s), \lfloor (r+p)/2 \rfloor, r$).
 \item[2.] If $f({\bf t}; s) * g_{[p,r]}(s)$ is  identically zero
then \\
{\tt findDegree}($f({\bf t}; s), \lfloor p/2 \rfloor, p-1$).
\end{itemize}
The call {\tt findDegree}($f({\bf t}; s), 0 , K$) takes
at most ${\rm log}(K)$ steps to find the degree.
Zero testing is done by substituting a positive vector
$(z_1, \ldots, z_d; w)$ where $0 < z_i \ll 1$ and $0 < w \ll 1$ which is
not a pole of $f({\bf t}; s) * g_{[p,r]}(s)$. Note that this Hadamard
product is a polynomial in $s$ whose coefficients have expansions of the
form $\sum_{b \in \Gamma} {\bf t}^b$ for some set $\Gamma$.
This means that the substitution of the above positive vector
gives a positive number unless the Hadamard product is identically zero.
\end{proof}

\noindent
{\sl Proof of Theorem \ref{complexity}: }
 By Theorem \ref{genfct} we can compute  the rational function
$G({\bf t}; s)$  in polynomial time. The degree of
$G({\bf t}; s)$ viewed as a polynomial in $s$
is ${\rm gap}(A,c)$. Theorem 17.2 of \cite{Sch} implies
that ${\rm gap}(A,c) \leq n D(A) \sum_{i=1}^n |c_i|$.
Let $K$ be this upper bound.
We observe that ${\rm log}(K)$ is a polynomial in the
bit size of the data. We can hence use
Lemma \ref{finddegree} to find ${\rm gap}(A,c)$.
\qed

\section{An Application to Algebraic Statistics}

We present an application to the
statistical theory of {\em disclosure limitation}.
See  \cite{CG} and \cite{DF} and the references  therein.
Suppose we are given data in the form of an
$n$-dimensional table of nonnegative integers. The aim is
to release some marginals of the table but not the
table's entries themselves. If the range of possible values that a
particular entry can attain in any table satisfying the released
marginals is too narrow, or even worse, consists of the
unique value of that entry in the actual table, then this entry may be
exposed. This shows the importance of determining tight integer upper
and lower bounds for each entry in a given table. A choice
of marginals corresponds to a simplicial complex on
$\{1,2,\ldots,n\}$ and can be represented by
a  zero-one matrix $A$, as described in
 \cite[\S 1]{HS}. In statistical language, the matrix
$A\,$ specifies a {\em hierarchical model} for a 
{\em contingency table} with $n$ {\em factors}.

The {\em table entry security problem} can be formally stated
as follows:
 suppose $u$ is a table with nonnegative integer entries,  where
  the marginals are computed according to a fixed hierarchical model $A$
  and let $u_{i_1 i_2 \cdots i_n}$
be a particular cell for such tables.  Compute
optimal lower and upper bounds $L$ and $U$ such that $L
  \leq u_{i_1 i_2 \cdots i_n} \leq U$ for all tables with the same 
marginals as $u$.

This problem is an integer program  (\ref{ip}):
minimize (or maximize)  $u_{i_1 i_2 \cdots i_n}$
 over all tables with nonnegative {\em
  integer} entries subject to fixing the marginals.
In view of the difficulty of solving integer programs
exactly, various researchers resorted to solving the
linear programming relaxation (\ref{lp}) instead:
 minimize (or maximize)  $u_{i_1 i_2 \cdots i_n}$ over all
tables with nonnegative {\em real} entries subject to fixing the
marginals. This relaxation is tractable,  but it usually fails
 to deliver the exact integers
$L$ and $U$.  One faces the problem of estimating the 
integer programming gap for the table security problem. 

We give the precise definition of the relevant matrices $A$.
Consider $d_1 \times \cdots \times d_n$-tables with 
entries $u_{i_1 i_2 \cdots i_n}$ where $1 \leq i_j \leq d_j$.
We fix a  hierarchical model by specifying  a collection
 of subsets $\,F_1, \ldots, F_k\,$ of $\{1, \ldots, n\}$.
The marginals of our table
are computed with respect to these subsets.
If $F_i = \{j_1, \ldots, j_s\}$ then the
$F_i$-marginal is a $d_{j_1} \times \cdots \times d_{j_s}$ table 
$b$ with  entries 
\begin{equation}
\label{umumum}
   b_{k_1 \cdots  k_s} \quad = 
    \sum_{i_{j_1} = k_1, \ldots, i_{j_s} = k_s} u_{i_1 \cdots i_n}. 
\end{equation}

\begin{ex}
 The classical {\em transportation problem} corresponds
to $d_1 \times d_2$-tables where the marginals are computed
with respect to $F_1 = \{1\}$ and $F_2 = \{2\}$. The 
{\em three-dimensional} transportation problem 
concerns $d_1 \times d_2 \times d_3$-tables with $F_1 = \{1,2\}$, 
$F_2 = \{1,3\}$, and $F_3 = \{2,3\}$. The marginals are
$$ b_{ij} = \sum_k u_{ijk}, \,\,\,  b_{ik} = \sum_j u_{ijk}, \,\,\,
b_{jk} = \sum_i u_{ijk}. $$      
For a discussion from the Gr\"obner basis perspective see
\cite[\S 14.C]{GB+CP}. \qed
\end{ex}

We define $A$ to be the zero-one matrix with
$d_1 d_2 \cdots d_n$ columns
that corresponds to the linear map that computes 
the marginals of tables. 
We let $u$ be  the vector of variables representing the cell entries.
Then $A \cdot u$ is the vector of the $k$ lower-dimensional
tables computed as in (\ref{umumum}).  There is a
transitive symmetry group acting on the columns of $A$, so
it suffices to examine the particular cell entry $u_{1 1 \cdots 1} $
which corresponds to the first column of $A$. 
The table entry security problem is the pair of integer programs
\begin{equation}
\label{newip}
{\rm Minimize} \,\, ({\rm Maximize})
  \,\,\,\, u_{11 \cdots  1} \quad
\hbox{subject to} \quad A \cdot u  =  b, 
\, u \geq 0, \, u \,\, \hbox{integral}.
\end{equation}

Our Theorem  \ref{firstthm} gives an exact formula
for the integer programming gap of these problems,
which we abbreviate by ${\rm gap}_-(A)$ and
${\rm gap}_+(A) $ respectively. Thus ${\rm gap}_+(A)$
is the worst error one gets when using
linear programming in computing
the bound $U$ for any $d_1 \times \cdots \times d_n$-table
with any fixed  margins $b$.

We illustrate our results for
$2 \times 2 \times 2 \times 2$-tables  $(u_{ijkl})$.
The {\em $K_4$-model} is specified by taking
all six two-dimensional  margins
$\, F_1 = \{1,2\}$,
$\,   F_2 = \{1,3\} $, 
$\,   F_3 = \{1,4\}$, 
$\,   F_4 = \{2,3\}$,  
$\,   F_5 = \{2,4\}$,  
$\,   F_6 = \{3,4\} $.
The zero-one matrix $A$ for the $K_4$-model
has $24$ rows and $16$ columns:
\renewcommand\arraystretch{0.6}
\fontsize{10}{12}
$$
A \quad = \quad \left[ \begin{array}{cccccccccccccccc}
 1 &  1 &  1 &  1 &  0 &  0 &  0 &  0 &  0 &  0 &  0 &  0 &  0 &  0 &  0 &  0 \\
 0 &  0 &  0 &  0 &  1 &  1 &  1 &  1 &  0 &  0 &  0 &  0 &  0 &  0 &  0 &  0 \\  
 0 &  0 &  0 &  0 &  0 &  0 &  0 &  0 &  1 &  1 &  1 &  1 &  0 &  0 &  0 &  0 \\
 0 &  0 &  0 &  0 &  0 &  0 &  0 &  0 &  0 &  0 &  0 &  0 &  1 &  1 &  1 &  1 \\
 1 &  1 &  0 &  0 &  1 &  1 &  0 &  0 &  0 &  0 &  0 &  0 &  0 &  0 &  0 &  0 \\ 
 0 &  0 &  1 &  1 &  0 &  0 &  1 &  1 &  0 &  0 &  0 &  0 &  0 &  0 &  0 &  0 \\
 0 &  0 &  0 &  0 &  0 &  0 &  0 &  0 &  1 &  1 &  0 &  0 &  1 &  1 &  0 &  0 \\
 0 &  0 &  0 &  0 &  0 &  0 &  0 &  0 &  0 &  0 &  1 &  1 &  0 &  0 &  1 &  1 \\
 1 &  0 &  1 &  0 &  1 &  0 &  1 &  0 &  0 &  0 &  0 &  0 &  0 &  0 &  0 &  0 \\
 0 &  1 &  0 &  1 &  0 &  1 &  0 &  1 &  0 &  0 &  0 &  0 &  0 &  0 &  0 &  0 \\
 0 &  0 &  0 &  0 &  0 &  0 &  0 &  0 &  1 &  0 &  1 &  0 &  1 &  0 &  1 &  0 \\
 0 &  0 &  0 &  0 &  0 &  0 &  0 &  0 &  0 &  1 &  0 &  1 &  0 &  1 &  0 &  1 \\
 1 &  1 &  0 &  0 &  0 &  0 &  0 &  0 &  1 &  1 &  0 &  0 &  0 &  0 &  0 &  0 \\
 0 &  0 &  1 &  1 &  0 &  0 &  0 &  0 &  0 &  0 &  1 &  1 &  0 &  0 &  0 &  0 \\
 0 &  0 &  0 &  0 &  1 &  1 &  0 &  0 &  0 &  0 &  0 &  0 &  1 &  1 &  0 &  0 \\
 0 &  0 &  0 &  0 &  0 &  0 &  1 &  1 &  0 &  0 &  0 &  0 &  0 &  0 &  1 &  1 \\
 1 &  0 &  1 &  0 &  0 &  0 &  0 &  0 &  1 &  0 &  1 &  0 &  0 &  0 &  0 &  0 \\
 0 &  1 &  0 &  1 &  0 &  0 &  0 &  0 &  0 &  1 &  0 &  1 &  0 &  0 &  0 &  0 \\
 0 &  0 &  0 &  0 &  1 &  0 &  1 &  0 &  0 &  0 &  0 &  0 &  1 &  0 &  1 &  0 \\ 
 0 &  0 &  0 &  0 &  0 &  1 &  0 &  1 &  0 &  0 &  0 &  0 &  0 &  1 &  0 &  1 \\
 1 &  0 &  0 &  0 &  1 &  0 &  0 &  0 &  1 &  0 &  0 &  0 &  1 &  0 &  0 &  0 \\ 
 0 &  1 &  0 &  0 &  0 &  1 &  0 &  0 &  0 &  1 &  0 &  0 &  0 &  1 &  0 &  0 \\  
 0 &  0 &  1 &  0 &  0 &  0 &  1 &  0 &  0 &  0 &  1 &  0 &  0 &  0 &  1 &  0 \\
 0 &  0 &  0 &  1 &  0 &  0 &  0 &  1 &  0 &  0 &  0 &  1 &  0 &  0 &  0 &  1 
\end{array} \right]  $$
\normalsize
Here the  cell entries are ordered lexicographically,
$u_{1111}, u_{1112}, \ldots, u_{2222}$.

We have the following computational result.

\begin{prop}
Every hierarchical model for $2 \times 2 \times 2 \times 2 $-tables
satisfies
$$ {\rm gap}_+(A) \quad \leq \quad 5/3 \quad = \quad 1.666666..., $$
and this bound is attained for the $K_4$-model, that is,
for the matrix $A$ above.
\end{prop}

\begin{proof}
We show that ${\rm gap}_+(A) = 5/3$ for the $K_4$-model.
Similar computations show that ${\rm gap}_+(A) < 5/3$
for all other simplicial complexes on $\{1,2,3,4\}$.

The monomial ideal $M(A,c)$ where $A$ is the above matrix and 
$c = [-1, 0, \ldots, 0]$ has $61$ minimal generators. Two of them are 
$$x_{1112}^3x_{1221}x_{1222}x_{2121}x_{2122}x_{2211}x_{2212} 
 \quad \mbox{     and     } \quad
x_{1122}x_{1212}x_{1221}x_{2111}x_{2222}^2.$$
This ideal has $139$ irreducible components. One of these  
components is 
$$
I(u, \tau) \, = \, \langle x_{1112}^2, x_{1121}^2, x_{1122}, x_{1211}^2, x_{1212}, x_{1221}, x_{1222}, \quad \quad \quad \quad \quad \quad \quad \quad \quad$$
$$ \quad \quad \quad \quad \quad \quad \quad \quad \quad
x_{2111}^2, x_{2112}, x_{2121}, x_{2122}, x_{2211}, x_{2212}, x_{2221}, x_{2222}^2 \rangle. $$

\noindent Here the non-zero components of $u$ are 
$u_{1112} = u_{1121} = u_{1211} = u_{2111} = u_{2222} = 1.
$
The gap value of $I(u,\tau)$  is $5/3$. 
We also give the optimal solution $v$ for the linear program (\ref{lp}) 
with $b = A \cdot u$:  
$$v_{1112} = v_{1121} = v_{1211} = v_{2111} = v_{2222} = 0,  $$
$$v_{1122} = v_{1212} = v_{1221} = v_{1222} = v_{2112} = 1/3, $$
$$v_{2121} = v_{2122} = v_{2211} = v_{2212} = v_{2221} = 1/3, \,\, 
v_{1111} = 5/3.$$
\end{proof}

The following corollary was pointed out to us by Rekha Thomas. 
\begin{cor} Let $M(A,c)$ be the monomial ideal corresponding 
to the minimization problem (\ref{newip}). Then
${\rm gap}_-(A,c) + 1$ is at least the maximum degree of $x_{11 \cdots 1}$ 
in any minimal generator of $M(A,c)$.  
\end{cor}
\begin{proof} This statement is equivalent to  
$${\rm gap}_-(A) \, \geq \, \max \{u_{11 \cdots 1}\,\, : I(u,\tau) \mbox{ irreducible component of } M(A,c)\}.$$
Let $I(u,\tau)$ be an irreducible component
such that $u_{11 \cdots 1} \geq 1$. Since $c = [1, 0, \ldots, 0]$, 
the objective function of the program (\ref{newgap}) corresponding to 
this component would be $b_{11 \cdots 1} \cdot t$. Moreover, by 
Theorem 2.5 in \cite{HT2}, the inequality 
$b_{11 \cdots 1} \cdot t \leq u_{11 \cdots 1}$ is a facet of the 
feasible region of (\ref{newgap}). Hence the optimal solution 
value is $u_{11 \cdots 1}$. Now Theorem \ref{firstthm} gives the result.  
\end{proof}

\begin{rmk} In the statistics literature there are 
various approaches to estimate $L$ and $U$ of the
table security problem. For general hierarchical models
an iterative algorithm for such an estimation is given in \cite{BG}.
A detailed analysis for {\em decomposable models} 
is given by Dobra and Fienberg \cite{DF}. 
\end{rmk}

\section{The gap fan}

In this section we allow the cost function $c$ to vary
in the programs (\ref{ip}) or (\ref{latprog}).
For each fixed matrix $A$ (resp.~fixed lattice ${\cal L}$)
we thus get a function
\begin{equation}
\label{gapfunction}
 {\rm gap}_A\, : \, \rr^n \longrightarrow \rr \,, \,\,
c \mapsto {\rm gap}(A,c).
\end{equation}
Our goal is to show that this function is piecewise linear, and
the natural pieces on which the function is linear form a fan
which we  call the {\em gap fan}.

Consider the {\em Gr\"obner fan} of
the matrix $A$. Following \cite{GB+CP} and \cite{ST},
this is the coarsest polyhedral fan in $\rr^n$ on which the
map $\, c \, \mapsto \, M(A,c)\,$ is constant.
Efficient software packages for computing
the Gr\"obner fan are described in \cite{HubT} and
\cite{CaTS}. The gap fan will be a refinement of
the Gr\"obner fan. Hence it suffices to describe
the gap fan on each Gr\"obner cone ${\mathcal K}$ separately.

We fix a maximal cone ${\mathcal K}$ of the Gr\"obner fan.
By the results of \cite{ST}, the
polyhedral cone ${\mathcal K}$ consists of
cost vectors such that the optimal solutions of all
integer programs (\ref{ip}) are constant as the right-hand-side vector
$b$ varies.
There exists a monomial ideal $M$
in $\rr[x_1,\ldots,x_n]$ such that
$$ {\cal K} \quad = \quad \bigl\{ \, c \in \rr^n \,\, : \,\,
M(A,c) = M \, \bigr\}. $$
Let  $\bigl\{I(u,\tau) \bigr\}$ be the finite set
of all  irreducible components of the monomial ideal $M$.
For each such component let $\,v_{u,\tau}\,$
be the optimal solution to the linear program
(\ref{gapprog}). It is not hard to see that these optimal
solutions are also constant as $c$ varies in ${\mathcal K}$.
Now we define the following polyhedral cone:
\begin{equation}
\label{gapcone}
{\mathcal G} \quad := \quad
\bigl\{(c,t) \in \rr^{n+1} \, : \, c \in {\mathcal K}  \hbox{ and }
c \cdot (u - v_{u,\tau}) \leq t \,\,\, \forall \,\,  I(u,\tau) \bigr\}
\end{equation}

\begin{thm} \label{lowenv}
 The projection of the lower envelope of ${\mathcal G}$
onto ${\mathcal K}$ is a polyhedral subdivision of ${\mathcal K}$.
The gap function is linear on each face of this subdivision.
\end{thm}

\begin{proof}
The facets of the lower envelope of ${\mathcal G}$
are of the form $\{(c,t) \in \rr^{n+1} \, : \,
c \cdot (u - v_{u,\tau}) = t\}$ where $t$ is the optimal
value of the program (\ref{gapprog})
for all $c$ such that $(c,t)$ is on this facet.
For such pairs $(c,t)$ we have $ t = {\rm gap}(A,c)$,
by Theorem \ref{firstthm}. Hence the lower envelope
of ${\mathcal G}$ is the graph of the gap function (\ref{gapfunction})
 over ${\mathcal K}$.
Now, clearly, the projection of the  lower envelope onto ${\mathcal K}$ is
a polyhedral subdivision of ${\mathcal K}$, and by its construction,
the  gap function (\ref{gapfunction}) is linear on each cone in
 this subdivision of ${\mathcal K}$.
\end{proof}

We define the {\em gap fan} as the refinement
of the Gr\"obner fan of $A$, where each Gr\"obner cone
${\mathcal K}$ is subdivided as in  Theorem \ref{lowenv}.
Our discussion implies:

\begin{cor} The gap function (\ref{gapfunction}) is a piecewise linear
function on $\rr^n$. It is linear on the cones
of the gap fan, but generally not on the Gr\"obner fan.
\end{cor}

\begin{ex}  Let us return to the example in the Introduction
but now with varying cost function. Our problem is
to make change using  pennies, nickels, dimes and quarters,
where the number of coins is fixed. The Gr\"obner fan of
the corresponding $2 \times 4$-matrix $A$
 has seven maximal cones. These cones and the
corresponding ideals can be computed using {\tt TiGERS}
\cite{HubT} or {\tt CaTS} \cite{CaTS}.

The gap fan has eight cones. Exactly one cone of
the Gr\"obner fan is divided into two cones. It
is the one corresponding to the ideal
in (\ref{earliermono}).
This Gr\"obner cone is defined by the inequalities
$$3n + q \geq 4d \,\,\, \mbox{ and } \,\,\, 8d \leq 5p + 3q.$$
The hyperplane defined by
$$ 305p - 135n - 308d + 138q = 0 $$
splits this cone into two pieces. On the  positive side of the hyperplane
the ${\rm gap}(A,c)$ is given by the irreducible component
$\langle p^5, n^3 \rangle$,
and on the negative side it is given by $\langle n^6, d^4, q \rangle$.
We list the irreducible components
of all seven initial ideals and the winning irreducible component
for each of them:
\renewcommand\arraystretch{1.0}
$$\begin{array}{|c|c|}
\hline
M(A,c) & \mbox{  Winning component(s)  } \\
\hline
\langle p^5, d^4 \rangle &  \langle p^5, d^4 \rangle \\
\hline
\langle p^5, q \rangle  \cap \langle p^5, n^3 \rangle \cap
\langle d^4, q \rangle &  \langle p^5, n^3 \rangle \\
\hline
\langle p^5, n^3 \rangle \cap \langle n^9, q \rangle &
\langle p^5, n^3 \rangle \\
\hline
\langle p^5, n^3 \rangle \cap \langle n^6, q \rangle  \cap
\langle n^3, q^2 \rangle &  \langle p^5, n^3 \rangle \\
\hline
\langle p^5, n^3 \rangle \cap \langle n^6, d^4, q \rangle  \cap
\langle n^3, q^3 \rangle &  \langle p^5, n^3 \rangle \\
                         &  \langle n^6, d^4, q \rangle \\
\hline
\langle n^6, d^4, q \rangle  \cap \langle n^3, d^8 \rangle &
\langle n^6, d^4, q \rangle \\
\hline
\langle n^6, d^4 \rangle &  \langle n^6, d^4 \rangle \\
\hline
\end{array}$$
\end{ex}

\end{document}